\documentclass[twocolumn,sloppy]{article}
\usepackage{graphicx}

\usepackage{latexsym}
\usepackage{amsfonts}
\usepackage{amsthm}
\usepackage{amsmath}
\def\L{\mathbb{L}}

\def\C{{\mathbb C}}

\def\Q{{\mathbb Q}}
\def\F{{\mathbb F}}
\def\A{{\mathbb A}}
\def\schk{{\operatorname{Sch}_k}}
\def\MotkQbar{{\operatorname{Mot_{k,\bar\Q}}}}
\def\Qbar{{\bar\Q}}
\def\frakg{{\mathfrak g}}
\def\frakh{{\mathfrak h}}
\def\fraka{{\mathfrak a}}
\def\strip{{\operatorname{strip}}}
\def\sgn{{\operatorname{sign}}}
\def\Ad{{\operatorname{Ad}}}
\def\ord{{\operatorname{ord}}}
\def\ac{{\operatorname{ac}}}
\def\id{{\operatorname{id}}}

\newtheorem{theorem}{Theorem}[section]

\newtheorem{conjecture}[theorem]{Conjecture}

\newtheorem{thesis}[theorem]{Thesis}

\newtheorem{proposition}[theorem]{Proposition}
\newtheorem{principle}[theorem]{Principle}
\newtheorem{problem}[theorem]{Problem}

\theoremstyle{definition}

\newtheorem{remark}[theorem]{Remark}

\begin{document}

\title{Can $p$-adic integrals be computed?}
\author{Thomas C. Hales}

\maketitle

\begin{abstract}
    This article gives an introduction to arithmetic motivic integration in
    the context of $p$-adic integrals that arise in
    representation theory.  A special case of the
    fundamental lemma is interpreted as an identity of Chow motives.
\end{abstract}

\section{Introduction}

This article raises a question in its title, and the short answer
to the question is that it still has not been
answered.\footnote{This article is based on a lecture at IAS,
    April 6, 2001
    http://www.math.ias.edu/amf/}\footnote{I would like to thank
    Carol  Olczak  for her assistance in preparing this
    manuscript.}%
    \footnote{I grant this paper to the public domain. No rights are reserved
        by the author.}
However, tools have now been developed to answer
questions such as this, and this article gives an introduction to some of these tools.

This article will concentrate on a particular family of integrals
that arise in connection with the representation theory of
reductive groups.  These are orbital integrals.  The clear
expectation is that these integrals can be computed, for reasons
that will be explained below.

This article will also touch on the fundamental lemma,
which is a conjectural identity that holds between certain orbital integrals.
This article will include a statement of the fundamental lemma in a special case.

The first sections may seem misplaced because they describes some methods that
are not in current use in representation theory, but by the end of the article, their
relevance will be established.

The central question in my research for some time is the question
of how to use a computer to prove theorems, particularly theorems
in geometry.  I hope to show that there is some interesting
geometry that arises in connection with $p$-adic integration, and
that computers can enhance our understanding of that geometry.

In Sections \ref{sec:tarski}, \ref{sec:integration}, and \ref{sec:motives}, three major
threads will be introduced: Tarski's decision procedure for the real numbers,
$p$-adic integration, and motives.  The other sections will tie these threads together
in the context of the fundamental lemma and $p$-adic orbital integrals.

\section{Tarski's decision procedure}
\label{sec:tarski}

Around 1930, Tarski proved a decision procedure for sentences in
the elementary theory of real closed fields.\footnote{An excellent
introduction to this topic, including a reprint of Tarski's
original article can be found in \cite{QECAD}.  A survey of recent
improvements in algorithms can be found at \cite{CRC} and
\cite{BPR}.}

Tarski's result can be formulated precisely in terms of a
first-order language. The language is built from the fifteen
symbols.

$$
    \begin{matrix}
    0 &1 &+ &*\\
    ( &) &= &<\\
    \forall &\exists &x &'\\
    \wedge &\vee &\neg\\
    \end{matrix}\\
$$
We will not go into the details of the syntax of the language.
\footnote{The general syntactic conventions of first-order languages can be found
    in \cite{enderton}.  A treatment of syntax in the context of algebraic
    structures can be found in \cite{fieldarithmetic}}
Each $x$ is followed by zero or more primes, and primes only occur
after $x$ or another prime.  We abbreviate $x$ followed by $n$
primes to $x_n$.  The language contains variables $x_n$, and the
constants $0$, $1$.  The variables and constants can be added and
multiplied (symbols $+$ and $*$).  Polynomial expressions can be
compared with the predicates $=$ and $<$. The quantifiers
($\forall$ and  $\exists$) should be understood as ranging over
the real numbers (or a complete ordered field).  For example, the
assertion that a quadratic polynomial has a root can be written in
this formal language as
\begin{equation}
    \neg (x'=0) \wedge \exists x (x'*x*x+x''*x+x'''=0)
    \label{eqn:quadratic}
\end{equation}
The formal language quickly becomes cumbersome, and we allow
ourselves certain informal shorthand conventions, for example,
writing Formula \ref{eqn:quadratic} as
$$
    a\ne 0 \wedge \exists x (a x^2 + bx +c=0),
$$
whenever a translation back into a formal statement of the language is clear.

Many things are noticeable absent from this little first-order
language.  There is no way to express particular real numbers in
this language such as $\pi=3.14159\ldots$, $e=2.71828\ldots$,
$\ln(2)$. There is no notion of set. There are no quantifiers that
range over subsets of the real numbers (for example, there are no
quantifiers over the integers).  There are no transcendental
functions such as the cosine function.  There is no calculus or
integration (except for formal derivatives of polynomials and the
like).

Tarski's result can be expressed as an algorithm for the
elimination of quantifiers in this first-order language. It takes
a formula in this language and manipulates it by an entirely
mechanical procedure into an equivalent form that contains no
quantifiers ($\exists$ $\forall$).  The formula that this
procedure gives as output is equivalent to the input in the sense
that the same $n$-tuples of real numbers satisfy the two formulas.

For example, if we apply Tarski's procedure to Formula
\ref{eqn:quadratic}, it returns something equivalent to the
quantifier-free formula
    $$
    \begin{array}{lll}
    &\neg  (x'=0) \wedge \\
    &(x''*x''-(1+1+1+1)*x'*x'' >0 \\
    &\vee x''*x''-(1+1+1+1)*x'*x'=0))
    \end{array}
    $$
or less formally,
    $$
    a\ne0 \wedge (b^2-4ac\ge0).
    $$
In other words, Tarski's procedure determines that a quadratic
equation has a real root if and only if the discriminant is
non-negative. In a similar way, the truth value of all sentences
in this language can be decided:  the truth value of an equivalent
sentence without quantifiers is trivially determined.

Here is a more difficult example, drawn from \cite[page 7]{QECAD}.
When is a quartic polynomial semi-definite?  Tarski's algorithm
takes the (formal translation of)
    $$\forall x (x^4 + p x^2 + q x + r \ge 0)$$
and returns a formula equivalent to
    $$
    \begin{array}{lll}
    (256 r^3 - 128 p^2 r^2 &+ 144 p q^2 r \\
                           &+ 16 p^4 r - 27 q^4 - 4 p^3 q^2 \ge0 \\
                     \wedge&\\
     \quad 8 p r - 9 q^2 - 2 p^3 &\le 0) \\
     \vee&\\
    \quad (27 q^2 + 8 p^3 \ge0
        \wedge &8 p r - 9 q^2 - 2p^3\ge 0)\\
         \wedge &\\
         \quad r\ge0&
    \end{array}
    $$

Tarski's original algorithm is very slow, but in 1975 George
Collins found a vastly improved method of quantifier elimination.
Further improvements are mentioned in the survey article
\cite{QECAD}.

The methods have improved to the point that the algorithms are of
practical importance. For instance, in robotics, quantifier
elimination can be used to determine whether two moving objects
will collide \cite{QECAD}.  Mathematica 4.0 implements an
experimental package in quantifier elimination \cite{mathematica}.
 There are highly nontrivial problems
in discrete geometry that can be expressed in this little
first-order language (for example, the dodecahedral conjecture
\cite{alg}).  The strategy is to squeeze non-trivial assertions
into this little language, and then let the general algorithms
prove the results.

\section{Pas's language}
\label{sec:pas}

This article is concerned, however, with $p$-adic quantifier elimination and not
with Tarski's quantifier elimination over the reals.  The first early results on
quantifier elimination can be found in articles by Ax-Kochen and Ershov
(\cite{AK} and \cite{E}).  The approach that we follow grows out of the article
{\it Decision procedures for real and $p$-adic fields} by Paul J. Cohen in 1969
(see \cite{cohen}).  Cohen's work on $p$-adic quantifier elimination was
refined and extended by various people (Denef \cite{denef}, Macintyre \cite{macintyre},
and Pas \cite{pas}).  We will describe $p$-adic quantifier elimination as it
is developed by Pas.

Pas defines a first-order language for complete Henselian rings
that is analogous to Tarski's first-order language for the theory
of complete ordered fields.   It contains the following tokens
    $$
    \begin{matrix}
    0 &1 &+ &*\\
    ( &) &= & <\\
    \forall &\exists &x &m &\xi &'\\
    \wedge &\vee &\neg\\
    \ord &\ac &\\
    \end{matrix}
    $$
The language consists of syntactically well-formed formulas in
this language. There are three sorts of variables $x$, $x'$, $x''$
(which we abbreviate to $x_0$, $x_1$, etc.), $m$, $m'$, $m''$
(which we abbreviate to $m_i$) and $\xi$, $\xi'$, $\xi''$, etc.
(which we abbreviate to $\xi_i$).

In the interpretations of this language, there are three algebraic structures:
a valued field (such as a $p$-adic field), a value group (the target of the valuation,
    which will typically be the additive group of the integers),
and a residue field. The variables $x_i$ are of the valued-field
sort, the variables $m_i$ are of the additive group sort, and the
variables $\xi_i$ are of the residue field sort. Correspondingly,
there are three sorts of quantification depending on the sort of
variable the quantified is attached to. The constant $0$ comes in
three sorts:  ($0_x$, $0_m$, and $0_\xi$). These are interpreted
as the zero element in the valued field, the additive value group,
and the residue field, respectively.  The addition symbol $+$ is
overloaded in that it is interpreted as addition in the valued
field, addition in the value group, or addition in the residue
field, according to its arguments.  (The syntax requires the
arguments to $+$ to be of the same sort.)

The function name {\it ord\/} is interpreted as the valuation on
the field.  If the model is a $p$-adic field, {\it ord\/} is
interpreted as the normalized valuation on the field.  The
function name {\it ac\/} is interpreted as an angular component
function.  On the units in the ring of integers, the
interpretation is the mapping from the units to its nonzero
residue in the residue field.  On general nonzero elements, it is
interpreted as the function that scales its argument by a power of
a uniformizer to make it a unit and then takes its image in the
residue field.  (Although a uniformizer is used to construct the
interpretation of the function {\it ac\/}, the uniformizer itself
does not appear in Pas's language.)  Expressions involving ``$<$''
are restricted to the additive group sort.

One of the design requirements of this language is that it be
small enough for there to be a quantifier elimination procedure.
By results of G\"odel, this would not be possible if the language
were to encompass the full arithmetic theory of the integers
\cite{godel}. For this reason, the language is restricted to the
additive theory of the value group.  That is, integer products
such as $m*m'$ are prohibited in the language.  Integer
expressions may be compared through equality and inequality ($=$
and $<$).   According to a result proved by Presburger in 1929, a
decision procedure exists for the additive theory of the integers
(\cite{pressburger}).

Just as in the case of the first-order theory of the reals, much
is missing from the language.  For instance, there is no
uniformizer in the language, so we cannot express $p$-adic
expansions of numbers in the valued field. As in the case of the
reals, there is no notation that would allow us to express sets in
this language.  It is impossible to express field extensions
directly (only indirectly through polynomials defining the roots,
for instance).  Most of Galois theory and local class field theory
will be inexpressible.

However, this language is small enough for there to be a procedure
of quantifier elimination.  In 1989, Pas, building on earlier
results, proved that the quantifiers of the valued field sort can
be eliminated, in the sense that an algorithm exists to produce an
equivalent formula without quantifiers of the valued-field
sort.\footnote{Pas's language gives quantifier elimination of
quantifiers of the valued field sort.  To eliminate all
quantifiers, Pas's result must be combined with Presburger's
quantifier elimination on the additive theory of the integers, and
with the theory of Galois stratification for quantifiers of the
residue field sort.} (Equivalence here means in the sense that for
any complete henselian ring with a residue field of characteristic
zero, the two formulas have the same set of
solutions.\footnote{Although Pas's procedure requires the residue
field to have characteristic zero,
    Pas, Denef, and Loeser are able to apply these results to $p$-adic fields with
    residue fields of positive characteristic.  This involves the use of
    ultrafilters and ultraproducts.  Finitely many primes are
    discarded in the process.})

One of the main applications of Pas's language and its quantifier
elimination procedure has been to the theory of $p$-adic
integration.  For example, Pas's original article contains results
about the Igusa local zeta function, which is a $p$-adic integral
(\cite{pas}).

\section{$p$-adic integration}
\label{sec:integration}

Let $F$ be a $p$-adic field of characteristic zero.  Let $\frakg$ be a
reductive Lie algebra defined over $F$, $X$ a regular semisimple element
of $\frakg(F)$.  Let $f$ be a function of compact support on $\frakg(F)$.
We consider the stable orbit $O^{st}(X)$ of $X$ (meaning the $F$-points of
the orbit of $X$ over an algebraic closure).  We pick an invariant measure $\mu$
on the orbit.  The integral of $f$ over $O^{st}(X)$ is called an orbital integral.
A fundamental problem is to compute
    $$\int_{O^{st}(X)} f\, d\mu.$$

These integrals arise repeatedly in the representation theory of
$p$-adic groups, in places such as the trace formula.  The
conjectural fundamental lemma (it will be discussed in Section
\ref{sec:FL}) is an identity of orbital integrals.  The fact that
the fundamental lemma has resisted all efforts to prove it is
closely related to the difficulty of computing orbital integrals.

\subsection{An example in $so(5)$}
\label{sec:so5}

An example, will illustrate the nature of these integrals.   Let
$\F_q$ be the residue field. Assume that its characteristic is not
$2$.   Let $\frakg=so(5)$. Assume that $X$ has that property that
the valuation of $\alpha(X)$ is independent of the root $\alpha$.
Assume that
    $$|\alpha(X)| =  q^{-r/2},$$
for an odd integer $r$.
Viewing $X$ as a linear transformation on a $5$-dimensional vector space, the roots
of the characteristic polynomial of $X$ are
    $$0, \pm t_1, \pm t_2.$$
Let $R_X$ be the quadratic polynomial in $k[\lambda]$ with roots
the reduction mod a uniformizer $\varpi_F$ of
    $$t^2_i/\varpi_F^r.$$
We have an elliptic curve $E_X$ over the finite field $k$ given by
    $$y^2 = R_X(\lambda^2).$$

There are test functions $f$ so that (for appropriate normalizations of measures)
we have
    \begin{equation}
    \int_{O^{st}(X)} f \,d\mu = A(q) + B(q) |E_X(k)|,
    \label{eqn:elliptic}
    \end{equation}
for some rational functions of $q$: $A$ and $B\ne0$. (See
\cite{hyperelliptic}).  The rational  functions $A$ and $B$ depend
on $f$. This special case gives an indication of what orbital
integrals can give.

What does it mean to calculate the orbital integral?  The naive
and completely unsatisfactory answer is that a calculation of an
orbital integral is to take a particular $p$-adic field, a
particular element $X$ and to program a computer to find the
complex number expressed on the right-hand side of Equation
\ref{eqn:elliptic}. A satisfactory answer to what it should mean
to calculate the orbital integral is to find the rational
functions $A$ and $B$, and to give the elliptic curve $E_X$. In
other words, what is really needed is a symbolic computation that
gets at the underlying variety (in this case an elliptic curve).
This is the sense in which I intend the question asked in the
title ``Can $p$-adic integrals be computed?''

If we examine this example more closely, we might ask what
features of the problem made this calculation possible?  The first
obvious feature is that as we vary the parameter $X$, the elliptic
curves do not change erratically; rather, they vary within a nice
family of elliptic curves over the finite field.

The second noteworthy feature is that as we move from local field to local field,
we obtain (in some sense) the ``same family'' of elliptic curves in each case.
It is this consistency as we go from one local field to another that makes it
reasonable to hope that a computer algorithm might be found to compute the orbital
integrals for all local fields.  We can view this as a single elliptic curve $E$ that
is defined over $\Q(a,b)$:
    $$y^2 = x^4 + a x^2 + b,$$
or as a family parameterized by $a$ and $b$.
All elliptic curves $E_X$ for all $p$-adic orbital integrals for all local fields
come as various specializations of this family.

To carry this example farther, we might look at some of the identities that are
predicted by Langlands's principle of functoriality.  One such identity (that is
needed for applications of the trace formula) predicts
an equality of orbital integrals between $so(5)$ and $sp(4)$.
It has the form
    $$
    \int_{O^{st}(X),so(5)} f\,d\mu = \int_{O^{st}(Y),sp(4)}f'\,d\mu'.
    $$
The data for $sp(4)$ is similar to that data for $so(5)$.  The elements $X$ and $Y$
are related through their characteristic polynomials $P_X$ (resp. $P_Y$):
    $$P_X(\lambda) = \lambda P_Y(\lambda).$$
(Note that $0$ is always a root of $P_X$.)
(The function $f'$ has to be related to $f$ in a suitable way.)
When these integrals are computed we find elliptic curves for $sp(4)$ as well, and
the identity of orbital integrals holds if and only if we have an identity of the
following form
    $$A(q) + B(q)|E_X(k)| = A(q)+B(q)|E'_Y(k)|.$$
It turns out that the elliptic curves $E_X$ and $E'_Y$ are not isomorphic (they
have different $j$-invariants), but they can be proved to have the same number of
points by producing an isogeny between $E_X$ and $E'_Y$.

This isogeny can be expressed as a single isogeny between two
elliptic curves $E$ and $E'$ over $\Q(a,b)$.  The conclusion of
this discussion is that this particular identity of orbital
integrals holds for all $p$-adic fields because of an identity of
two Chow motives over $\Q(a,b)$: that is, $E$ is isogenous to
$E'$.

What this suggests is a general hope that there are global objects attached to $p$-adic
integrals.  The global object should be something like a Chow motive.  Identities
of $p$-adic integrals should be consequences of identities of Chow motives.

\section{Motives}
\label{sec:motives}

We are finally in the position to give a precise definition of what it means to compute
a $p$-adic integral.  We state it as a thesis:

\begin{thesis} The computation of a $p$-adic integral is an effective algorithm to
obtain the underlying virtual Chow motive.
\end{thesis}

This thesis is incoherent unless virtual Chow motives are associated with general
families of $p$-adic integrals.  That this should be so was articulated by Loeser
in Strasbourg in 1999 \cite{strasbourg}.

\begin{principle} (Denef-Loeser Principle)  All ``natural'' $p$-adic integrals are
motivic.
\end{principle}

Without committing Denef and Loeser to any particular definition of ``natural,''
as representation theorists, we like to think that the important integrals that
arise in representation theory
are natural.  We are thus led to an investigation of motivic underpinnings of
$p$-adic integrals.

\subsection{An Example of Motivic Integration}
\label{sec:motivic}

Motivic integration was introduced by Kontsevich in a lecture in
Orsay in 1995 \cite{kontsevich}.  The properties of motivic
integration have been developed in a series of fundamental
articles by Denef and Loeser \cite{denefloeser1},
\cite{denefloeser2}, \cite{denefloeser3}, \cite{denefloeser4}. In
fact, my entire article is nothing but an application of the
beautiful circle of ideas that they develop.

To describe the theory in a few words, I will describe motivic integration by analogy
with $p$-adic integration.  Consider the following elementary $p$-adic integral:
    $$
    \begin{array}{lll}
    \int_{\F_q[[t]]} |x|^m\,dx &=
    \sum_{\ell=0}^\infty |\varpi^\ell|^{m+1}\int_{|u|=1}\frac{du}{|u|}\\
    &=(1 + q^{-(m+1)}+ q^{-2(m+1)}\cdots)(1-q^{-1})\\
    \end{array}
        $$
The answer is independent of the field, so we are tempted to write for any field
(say a field of characteristic zero):
    $$
    \begin{array}{lll}
    \int_{k[[t]]} |x|^m\,dx &=
    \sum_{\ell=0}^\infty |\varpi^\ell|^{m+1}\int_{|u|=1}\frac{du}{u}\\
    &=(1 + q^{-(m+1)}+ q^{-2(m+1)}\cdots)(1-q^{-1})\\
    \end{array}
    $$
The only difficulty is in the interpretation of $q$.  Kontsevich
supplies the answer with motivic integration: it is a symbol.
More specifically, in the case of finite fields
 it is a symbol attached to the affine line $\A^1$, and for general fields we can
continue to view it as a symbol attached to the affine line.
To mark the change of context from $p$-adic fields to more general fields, we replace
the symbol $q$ with $\L$ (to suggest the Lefschetz motive).

\subsection{Rings of virtual motives}

Section \ref{sec:motivic} describes a simple example of motivic
integration.  You integrate much in the same way as with $p$-adic
integration, but whenever in $p$-adic integration it becomes
necessary to count points on a variety, with motivic integration
you introduce a new symbol for that variety and move on. Although
a single symbol ($q$ or $\L$) suffices for the one example that
was shown, motivic integration in general will require a host of
new symbols. Integration should be linear, so relations must be
introduced among the symbols to make motivic integration linear.
In that example, $q$ (or $\L$) occurs as a denominator, and this
will require us to invert $\L$ in the ring we construct. In that
example, the answer is a limit (that is, an infinite sum), and
this will require us to complete the ring we construct.


\section{Rings of motives}

Let $k$ be a field of characteristic zero.
Let $\schk$ be the category of varieties over $k$.
Let $K_0(\schk)$ be the Grothendieck ring of varieties over $k$.
It is the commutative ring generated by symbols
    $$
    [S]
    $$
for each variety $S$ over $k$.   The relations are
    $$
    [S\times S']=[S][S']
    $$
and if $S'$ is closed in $S$, then
    $$
    [S] = [S\setminus S'] + [S'].
    $$

We let $\L =[\A^1]$.
Let $K_0(\schk)_{loc}$ be the ring obtained by inverting $\L$.

We let $\MotkQbar$ be the category of Chow motives over $k$
with coefficients in $\Qbar$, an algebraic closure of $\Q$.
(This category is described in detail in \cite{scholl}.)
The objects in this category are triples
    $$
    (S,p,n)
    $$
where $S$ is a variety over $k$, $p$ is a projection operator over
    $\Qbar$, and $n$ is an integer.
This category is an additive category, but it is not abelian
\cite{scholl}.

Let $K_0(\MotkQbar)$ be the Grothendieck group of the additive
category $\MotkQbar$.  The generators of this group is the set of
objects in $\MotkQbar$. By a fundamental result of Gillet and
Soulet \cite{soulet} and Guill\'en and Navarro Aznar \cite{aznar},
there is a  homomorphism of rings
    $$K_0(\schk)\to K_0(\MotkQbar)$$
that takes the symbol $[S]$ of a smooth projective variety $S$ to
    the generator associated with $(S,\id,0)$, where $\id$ is the
    identity projection operator.
The image of $\L$ under this homomorphism is invertible. Thus, the
homomorphism extends to $K_0(\schk)_{loc}$. Let
$K_0^v(\MotkQbar)_{loc}$ be the image of this homomorphism.

There is a filtration $F^m K_0^v(\MotkQbar)_{loc}$ on this group
given by $S/\L^i \in F^m$ iff $\dim S - i \le -m$.  We let $\hat
K_0^v(\MotkQbar)_{loc}$ be the completion of with respect to this
filtration.  This is the ring  in which motivic integrals take
their values (or sometimes its tensor product with $\Q$).

\subsection{Arithmetic motivic integration}

In 1999, Denef and Loeser developed an arithmetic theory of
motivic integration in \cite{DS}.  (This theory is distinct from a
geometric theory of motivic integration that was developed
earlier.) In their article, that Denef and Loeser describe the
three threads introduced in Sections \ref{sec:tarski},
\ref{sec:integration}, and \ref{sec:motives}, and show how they
relate to one another. In their article, they make two fundamental
discoveries:
    \begin{enumerate}
    \item Motives can be attached to formulas in Pas's language.
    \item The trace of Frobenius on the motive equals the
    $p$-adic integral over the $p$-adic set defined by the
    formula.
    \end{enumerate}

\begin{figure}[htb]
  \centering
  \includegraphics{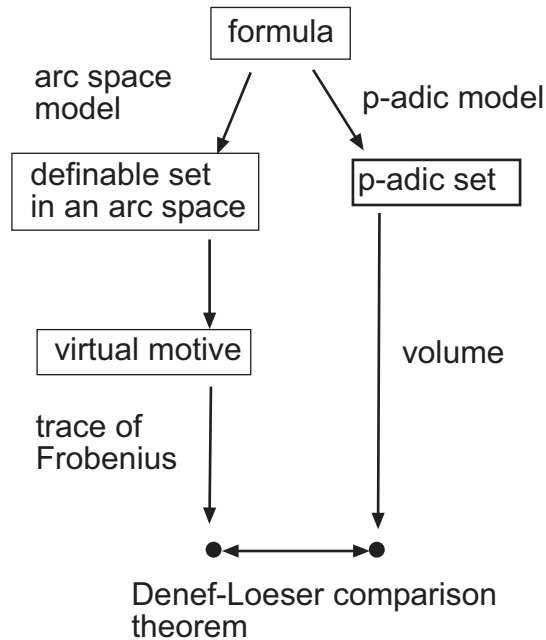}
  \caption{The Denef-Loeser comparison theorem}
  \label{fig:dlcomparison}
\end{figure}

The process is represented schematically in Figure
\ref{fig:dlcomparison}.
A formula in Pas's language can be interpreted two ways, leading
to two different integrals. First the formula can be interpreted
over a $p$-adic field.  The $p$-adic set of points that satisfy
the formula has a volume. The formula can also be interpreted over
a henselian field (such as $\C((t))$).  The set of points that
satisfy the formula has a motivic volume (an element of the ring
$\hat K_0^v(\MotkQbar)_{loc}$).

The comparison theorem of Denef and Loeser asserts that the trace
of Frobenius against this virtual motive is equal to the $p$-adic
volume of the $p$-adic set.  In particular, if a $p$-adic set has
the special form given by the set of points satisfying a formula
in Pas's language, then a motive can be attached to it. It is this
comparison theorem that will permit us to show that interesting
$p$-adic integrals have a motivic interpretation.

\section{The fundamental lemma}

All the sections until now have been an extended introduction to
provide context for the results I am about to describe.  Roughly
speaking, I have found that orbital integrals can be placed into
the framework of Denef and Loeser.

\subsection{Strips}
\label{sec:strips}

Let $F$ be a $p$-adic field of characteristic $0$.
Let $k$ be the residue field of $F$.
Fix parameters $n$, $k$, and $r$ satisfying the following conditions.
    \begin{itemize}
    \item $n$ is a positive integer.
    \item $k$ is an integer $k\le n$.
    \item $r$ is a rational number.  Write it as $r = \ell/h$,
    with $\ell$ and $h$ relatively prime.
    \end{itemize}

Let $\frakg =so(2n+1)$. There are endoscopic Lie algebras
    $$\frakh = so(2k+1)\times so(2n-2k+1).$$
That is, we take a product of two orthogonal Lie algebras, whose
ranks add up to that of $\frakg$.
(Endoscopy was originally defined in terms of groups, but it has
become common practice to follow the practice of Waldspurger and
to pass to the Lie algebras.)

We define a subset of $\frakg$ that I will call a strip.  It
depends on the parameter $r$.
Define $\strip(r)$ to be the set of all $X\in\frakg$ such that
$|\alpha(X)|=q^{-r}$ for all roots $\alpha$.  These elements
are called elements {\it equal valuation}.

Write the characteristic polynomial $P_X(\lambda)$ as
    $$P_X(\lambda) = \lambda P^0_X(\lambda).$$
Let $\bar \F_q$ be an algebraic closure of $\F_q$. Let $R_X$ be
the separable polynomial in $\F_q[\lambda]$ with roots in
    $\bar\F_q$
given by the reduction of the elements
    $$t_i^h/\omega_F^\ell,$$
where $t_i$ are the nonzero roots of $P_X$ (that is, the roots of
$P_X^0$).  The elements $t_i^h$ have been multiplied by an
appropriate power of the uniformizer, so that they become units.
As a result, the roots of $R_X$ are nonzero.

\subsection{Aside on Equal Valuation}
\label{sec:aside}

In this article, we do not justify our restriction to this special
kind of elements.  Without going into the details, it seems
to me that
the study of orbital integrals can be divided into two quite
different parts.  The part discussed in this article is that
of elements of equal valuation.  It seems that geometric methods
such as motivic integration are very important for this part.

For elements of nonequal valuation, it seems that a quite
different set of methods will be relevant.  Here issues such as
homogeneity (generalizing the results of Waldspurger
\cite{waldhom} and DeBacker \cite{debacker}) and descent for
orbital integrals \cite{descent} should be relevant.

This is currently merely speculation, but it is the reason that I
am restricting to elements of equal valuation.  It seems that
different groups could study these two different kinds of orbital
integrals with little interaction and few shared methods.

\subsection{A hope}

\begin{conjecture} If $X$ and $X'$ are elements in
$\strip(r)$ such that $R_X=R_{X'}$, then their orbital integrals
are equal. \label{conj:locallyconstant}
\end{conjecture}

\begin{figure}[htb]
  \centering
  \includegraphics{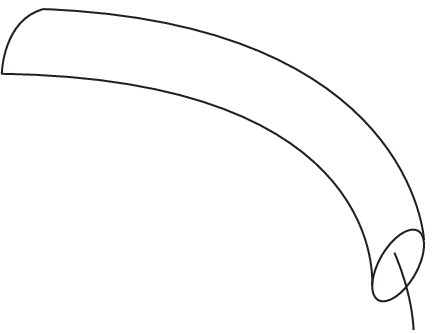}
  \caption{A tube}
  \label{fig:strip}
\end{figure}

We represent $\strip(r)$ schematically as a long rectangular strip
(a union of semisimple orbits).  Around the conjugacy class of $X$
we can draw a tube (a thickened neighborhood) of all the elements
in the strip with the same reduced characteristic polynomial
$R_X$.  (Figure \ref{fig:strip}.)
The function $X\mapsto R_X$ thus partitions $\strip(r)$ into tubes.

\subsection{A $p$-adic fundamental lemma}
\label{sec:FL}

Langlands and Shelstad define a transfer factor $\Delta(X,Y,Z)$
on $so(2n+1)\times so(2k+1)\times so(2(n-k)+1)$.  On the strip $r$
it has the form
    $$q^c \sgn(X,Y,Z)$$
for some constant $c=c(n,k,r)$ and some function $\sgn$ taking
values in
    $\{0,1,-1\}$.
Let $O_F$ be the ring of integers of the $p$-adic field $F$. They
conjecture that for appropriate normalizations of measures
\cite{hales-simple}, we have the following special case of the
fundamental lemma
\begin{conjecture} (Langlands-Shelstad)
For all $Y$ and $Z$ regular semisimple such that there exists a
regular semisimple $X$ in $\frakg$ such that $P_X^0=P_Y^0P_Z^0$,
we have
    $$
    \begin{array}{lll}
    q^c \sum_X \int_{O(X)\cap \frakg(O_F)}\sgn(X,Y,Z) &=\\
    \quad\int_{O^{st}(Y)\times O^{st}(Z)\cap \frakh(O_F)} 1.&\\
    \end{array}
    $$
\end{conjecture}

The sum runs over representatives of all regular semisimple
conjugacy classes. The function $\sgn(X,Y,Z)$ is zero unless
$P_X^0=P_Y^0 P_Z^0$. It is enough to restrict the sum to such
representatives.



\begin{proposition} If $F$ is a field of sufficiently large
residual characteristic, then the sign of the transfer factor in
$so(2n+1)$,
    $$\sgn^{-1}(x), \text{ for } x\in\{0,1,-1\},$$
is given by a formula in the language of rings.
\end{proposition}

The proof will be given in a separate article.  The surprising
thing about this calculation is that the full strength of Pas's
language is not required.  That is, the transfer factor is
expressed without the functions $\ord$ and $\ac$, and without
quantifiers over the additive group and residue field.  This means
that we can define a transfer factor for any field.

The formula for $\sgn^{-1}\{-1,1\}$ does not require quantifiers.
It is the set of elements $(X,Y,Z)$ with
    $$\lambda P_X = P_Y P_Z$$
for which $X$ is regular (which is expressed as the nonvanishing
of the resultant
    $$\operatorname{resultant}(P_X,P'_X)\ne0.$$

The starting point for the proof of the proposition is
Waldspurger's simplified formula for the transfer factors on the
Lie algebra of classical groups \cite{waldspurger}.

From this proposition, the Denef-Loeser construction gives us a
    $+1$-motive in $\hat K_0^v(\MotkQbar)_{loc}$.
It also gives a $-1$-motive in the same set.  In this sense, we
can affirm that the Langlands-Shelstad transfer factor is a
motive.

\subsection{Orbital integrals}

A serious difficulty that we encounter in the study of orbital
integrals is that individual orbits of semisimple elements are not
given by a formula in Pas's language.  In fact, the characteristic
polynomial
    $$P_X \in F[\lambda]$$
has $p$-adic coefficients, which cannot be expressed in the
language.  (Without a uniformizer in the language, we cannot
express the $p$-adic expansion of the coefficients.)

Our only hope is to use the fact that orbital integrals are
locally constant.  We place each orbit into a larger tube, where
the tube is large enough to be defined by a formula in Pas's
language.  Each tube is defined by the set of elements with a
given reduced characteristic polynomial with coefficients in the
residue field (see Section \ref{sec:strips}):
    $$R_X \in \F_q[\lambda].$$

To patch these together into a global object, we let $k$ be a
finite extension of $\Q$, with ring of integers $O_k$.   Each
polynomial $R_X$ is a specialization of a polynomial
    $$\dot R_X\in S[\lambda],$$
where $S$ is the coordinate ring over $O_k$ of the set of regular
orbits
    $$Z = Z_r = \fraka/\Ad\, A$$
on an appropriate Lie algebra $\fraka$ and group $A$, defined over
$O_k$. We can then use $\dot R_X$ to define a formula in
    $${\cal L}_{pas}(S).$$
(This denotes the Pas's language extended by a symbolic constant
for each element of $S$, as described in \cite[Section
6.3]{fieldarithmetic}.)

The formula for the set of elements in the tube with transfer
factor equal to $+1$ gives, by the construction of Denef and
Loeser, a Chow motive
    $$\Theta^{G,+}_{n,k,r}.$$  The negative part
of the tube gives a second Chow motive
    $$\Theta^{G,-}_{n,k,r}.$$
The tube on the endoscopic groups gives a third Chow motive
    $$\Theta^{H,st}_{n,k,r}.$$
Recall that the $p$-adic transfer factor has the form
    $$\pm q^{c}$$ for some constant $c=c(n,k,r)$.
We are thus able to formulate a motivic fundamental lemma:

\begin{conjecture}  Given $n$, $k$, and $r$, we have
    $$\L^c (\Theta^{G,+}_{n,k,r} - \Theta^{G,-}_{n,k,r})
        = \Theta^{H,st}_{n,k,r}$$
in
    $$\hat K_0^v(M_{\Q(Z_r),\bar \Q})_{loc,\Q}$$
\end{conjecture}

This single identity of Chow motives governs the fundamental lemma
over the entire $\strip(r)$ at almost all places.

\begin{remark} The Denef-Loeser comparison theorem relates the
trace of Frobenius on these motives to the traditional fundamental
lemma.  The Denef-Loeser comparison theorem in its current form is
not quite strong enough to deduce the fundamental lemma from its
motivic form.  However, I hope that these are relatively minor
obstacles that future research should be able to surmount.

First of all, we need the Denef-Loeser comparison theorem for the
finitely generated extension $\Q(Z_r)/\Q$.  Denef and Loeser give
two comparison theorems, one for $p$-adic integration on local
fields of positive characteristic, and another on local fields in
characteristic zero.  The comparison theorem in characteristic
zero assumes that the field is a finite extension of $\Q$ and
hence it cannot be applied to $\Q(Z_r)$.\footnote{Denef and Loeser
(private communication) have informed me that they can relax this
restriction for $p$-adic fields of characteristic zero that are
unramified over $\Q$.  This is expected to be sufficient for
applications to the fundamental lemma.}

The second restriction of the Denef-Loeser comparison theorem is
that if $R$ is a normal domain with field of fractions $\Q(Z_r)$,
then the comparison theorem holds at all closed points $x$ of
    $\operatorname{Spec}\,R_f$ for some non-explicit element $f$
of $R$.  It is possible that $f$ blocks a comparison of some
elements of the $p$-adic field.
\end{remark}

\begin{remark} One of the most interesting aspects of this
calculation is that it shows that there are two local-global
pathways.  The local-global pathway connecting automorphic
representation theory with the representation theory of local
fields is a well-established part of the Langlands program.  It
generalizes the pathway between global and local class field
theory.  The Denef-Loeser apparatus gives a genuinely new pathway
between global and local objects.  (To see that it is a different
pathway, observe that a global regular semisimple element in a
reductive group
    $$\gamma\in G(\Q)$$
lies in an unramified Cartan subgroup at almost every place.
However, the Denef-Loeser construction quite often globalizes
local data that is everywhere ramified.)

Thus, we have global problems in automorphic representation theory
that are localized by the first pathway, and then globalized again
by the second pathway to obtain conjectural identities of Chow
motives.
\end{remark}

\section{Open Problems}

We conclude this article with four problems that are raised by the
study of $p$-adic integrals from the vantage point of motivic
integration.

\begin{problem} Give effective algorithms to find the Chow motives
    $$\Theta^{*,*}_{n,k,r}.$$
\end{problem}

By solving this problem, we succeed in computing $p$-adic orbital
integrals in the sense proposed in this article.   The quantifier
elimination procedures (Pas's algorithm \cite{pas}, Presburger's
algorithm \cite{enderton}, and Galois stratification
\cite{fieldarithmetic}) are entirely algorithmic.  Thus, I hope
that this first problem can be settled.

\begin{problem} Prove the {\it hope} of Conjecture \ref{conj:locallyconstant}:
if $R_X = R_{X'}$, then
the orbital integrals of $O^{st}(X)$ and $O^{st}(X')$ are equal.
\end{problem}

In unpublished work, Clifton Cunningham has made progress toward a
solution of this second problem.

\begin{problem} Extend the results to degenerate elements $X$ that
do not lie in any strip $r$.  In particular, find {\it finitely}
many motives over finitely generated extensions of $\Q$ that
govern the fundamental lemma for all $X\in\frakg$ over almost all
completions of any number field.
\end{problem}

This seems to me to be a difficult problem.  As an earlier section
states (see \ref{sec:aside}), the methods involved here seem to be
methods of generalized homogeneity laws in the spirit of
Waldspurger and DeBacker (\cite{waldhom} and \cite{debacker}).

\begin{problem} Prove the motivic fundamental lemma.
\end{problem}

If the first three problems can be solved, then we have an
algorithm to compute the Chow motives that govern the fundamental
lemma for a given group.  However, more is needed.  First there is
the problem of equality: given two Chow motives, is there an
algorithm to determine if they are equal?

Second, there is the problem of induction.  Even if we have an
algorithm to check the fundamental lemma for one group, how do we
give a proof for all reductive groups at once?  Here it seems to
me that we need to develop a deeper understanding of the motives
that arise in connection with the fundamental lemma.

\section{Conclusion}

The Denef-Loeser apparatus of arithmetic motivic integration seems
to mesh well with certain $p$-adic integrals that arise in
representation theory.

We should investigate how far motives permeate representation
theory of $p$-adic groups.  If we believe with Denef and Loeser
that all natural $p$-adic integrals are motivic, then the
influence of the motivic point of view will be far-reaching.  One
can speculate that many of the basic objects of representation
theory (such as Harish-Chandra characters) have a motivic nature.

The hope is that the motivic interpretation will allow us to
calculate $p$-adic integrals that have resisted all efforts until
now.

    \renewcommand{\thefootnote}{}
    \footnote{version 5/19/02.}
    \footnote{Research supported in part by the NSF}

\end{document}